# Urban Freight Transportation Planning: A Dynamic Stackelberg Game-Theoretic Approach

**Bo Zhang, Tao Yao, Terry L. Friesz and Hongcheng Liu**
Department of Industrial & Manufacturing Engineering, Pennsylvania State University, University Park, PA 16802, United States

**Abstract** In this paper we propose a dynamic Stackelberg game-theoretic model for urban freight transportation planning which is able to characterize the interaction between freight and personal transportation in an urban area. The problem is formulated as a bi-level dynamic mathematical program with equilibrium constraints (MPEC) which belongs to a class of computationally challenging problems. The lower level is dynamic user equilibrium (DUE) with inhomogeneous traffic that characterizes traffic assignment of personal transportation given the schedule of freight transportation. The upper level is a system optimum (SO) freight transportation planning problem which aims at minimizing the total cost to a truck company. A mathematical program with complementarity constraints (MPCC) reformulation is derived and a projected gradient algorithm is designed to solve this computationally challenging problem. Numerical experiments are conducted to show that when planning freight transportation the background traffic is nonnegligible, even though the amount of trucks compared to other vehicles traveling on the same network is relatively small. What's more, in our proposed bi-level model for urban freight transportation planning, we find a dynamic case of a Braess-like Paradox which can provide managerial insights to a metropolitan planning organization (MPO) in increasing social welfare by restricting freight movement.

## 1. INTRODUCTION

Urban freight transportation, sometimes also referred to as city logistics, aims to reduce negative externalities, such as emission, noise and congestion, associated with freight activity while supporting their economic and social development (Crainic et al. 2009). Today the problem is even more important and challenging with the growing number of private vehicles, soaring demand of urban freight transportation services and increased recognition of the need for a paradigm shift toward environmentally sustainable logistics and freight technologies.

Urban freight transportation has attracted lots of research efforts in the past decades from different perspectives including transportation regulation, emission estimation and reduction, transportation planning, etc. Most studies about urban freight transportation planning modeled the problem based on a framework of the vehicle routing problem (VRP), with the objective of minimizing the total cost/delay to the carriers or truck companies while satisfying demand consumption constraints. These studies mainly focus on modeling freight activities. However, personal transportation and its impact on the freight transportation planning haven't been well studied.

To fill the gap mentioned above, in this paper we propose a dynamic game-theoretic model for urban freight transportation planning which is able to characterize the interaction between freight and personal transportation. Specifically, the problem is formulated as a dynamic mathematical program with equilibrium constraints (MPEC) which can handle the dynamic traffic assignment of inhomogeneous road users. Theoretical properties of this new model are discussed. To achieve tractability and guarantee solution quality, a mathematical program with complementarity constraints (MPCC) reformulation is derived. The problem finally can be solved by a proposed projected gradient algorithm. Numerical results show that the interaction of freight and personal transportation is nonnegligible to the truck company. Furthermore, by extensive numerical tests, we find a Braess-like paradox based on our proposed model which can provide managerial insights to a metropolitan planning organization (MPO) so that the social welfare may be improved by appropriate policy-making.

The paper is organized as follows. In Section 2, we present a comprehensive literature review and discuss how this paper distinguishes from and contributes to the existing body of literature. In Section 3, a dynamic Stackelberg game-theoretic model is built for urban freight transportation planning and some theoretical results are provided. In Section 4, we reformulate the problem as a mathematical program with complementarity constraints (MPCC) and design a projected gradient algorithm to solve it. Numerical experiments are presented and the results are analyzed in Section 5. In the end, Section 6 concludes the paper.

## 2. LITERATURE REVIEW

The intercity freight transportation assignment problems that focus on the forecast of transportation flows were intensively studied in 1970s and 1980s. Friesz et al. (1983) provided a survey of the predicative intercity freight network models and discussed the advantage of combined shipper-carrier models and spatial equilibrium models. Friesz et al. (1986) and Harker and Friesz's (1986a, 1986b) proposed to sequentially and simultaneously model and solve a shipper-carrier static freight assignment problem, respectively. More recently, Agrawal and Ziliaskopoulos (2006) proposed a dynamic shipper-carrier freight assignment model and an iterative approach to solve for the shippers' equilibrium solution. Nevertheless, due to the need in modeling each shipper and carrier's sub-network, the shipper-carrier approach is difficult to implemented practice (Crainic et al. 2007). Different from the shipper-career approach, some transshipment network assignment models were proposed to solely address carriers' transport operations (See Guélat et al. (1990) and Chow and Ritchie (2012) for some examples).

According to Crainic (2000), intercity freight transportation, also known as long-haul transportation, is a tactical (medium-term) planning issue for a transportation agency, while urban freight transportation (city logistics) is more of an operational (short-term) planning issue which is mostly modeled under the framework of a vehicle routing problem (VRP). The common objective of such models are to minimize the total cost/delay to the carriers or truck companies while satisfying demand consumption constraints by designing optimal pick-up or delivery routes from one or several nodes (depots) to a set of other nodes in a transportation network. Vehicle routing problem with time windows (VRPTW) is an extension of VRP which incorporates time constraints for freight pick-up and delivery. Some discussions or surveys of VRP/VRPTW can be found in, e.g. Kulkarni and Bhave (1985), Golden and Assad (1986), Laporte (1992) and Kallehauge et al. (2005). Specifically, in urban freight transportation planning studies, Taniguchi and Thompson (2002) studied a stochastic VRPTW which incorporated travel time variance. Crainic et al. (2009) proposed an integrated model that addresses short-term scheduling of operations and resource management based on a two-tiered distribution structure. Two new problem classes which are extensions of VRPTW were introduced and possible solution avenues were discussed. The two-tiered city logistics model was further extended to address demand uncertainty in Crainic et al. (2011). However, as Ambrosini and Routhier (2004) and Paglione (2006) noted, among the urban freight transportation planning literature, there is a lack of behavioral models that characterize the interactions of private economic and transport agents.

Nonetheless, the study of competition and cooperation among different road users is not new. Yang et al. (2007) incorporated the routing behaviors of system optimum (SO), user equilibrium (UE) and Cournot-Nash (CN) travelers using a static Stackelberg game with perfect information. In particular, the SO traveler is the leader and the UE and CN travelers are the followers. The UE and CN travelers make their routing decisions in a mixed equilibrium behavior given the SO traveler's routing decision, while the SO traveler optimizes its routes considering the potential reactions of UE and CN travelers towards its routing decision. Since urban freight transportation planning deals with short-term operations and planning issues, a model that can dynamically characterize interactions of road users is necessary.

Dynamic user equilibrium (DUE) captures the routing behaviors of individual travelers in a spatial network in a way that the effective unit delay/cost, including early/late arrival penalties, of traveling on all utilized path at any departure time is identical (see Friesz (2010) for a detailed discussion of DUE). As the first to model timely interactions of freight and personal transportation, this paper presents a dynamic Stackelberg game in which the leader is a truck company aiming at optimizing the freight transportation and the followers are individual travelers whose travel behaviors follow DUE with inhomogeneous traffic. The formulation belongs to a challenging set of mathematical programs which is known as dynamic mathematical program with equilibrium constraints (MPEC). Moreover the time shifts in the DUE model make solving this dynamic MPEC even more challenging. Friesz and Mookherjee (2006) proposed to use an implicit fixed point algorithm to accommodate the time shifts and based on this idea. Friesz et al. (2007) discussed two algorithms to solve a dynamic optimal toll problem with equilibrium constraints (DOTPEC) which is a specific type of dynamic MPEC. Yao et al. (2012) introduced toll price uncertainty into the DOTPEC and proposed a bi-level heuristic method to solve the resulting robust DOTPEC. Distinct from the DOTPEC, the lower level of our model is DUE with inhomogeneous traffic which requires a refinement of the traffic dynamics and network loading. Whence, in this paper, we propose a new dynamic MPEC model and discuss its theoretical properties. Moreover, to balance the quality of solutions and computational efficiency, instead of using heuristic algorithms we propose a MPCC reformulation and design a projected gradient algorithm to unlock the problem.

## 3. DYNAMIC STACKELBERG GAME-THEORETIC MODELING

In this study, we consider two types of traffic: freight transportation by trucks and personal transportation by private vehicles. Note that public transportation by buses is ignored in this model since its travel behavior (departure time and route choice) is almost fixed and independent of the behaviors of other vehicles. So in this paper, by personal transportation, we refer to private vehicles only. Trucks are controlled by a truck company who aims at minimizing the total transportation cost/delay while satisfying the travel demand (i.e., the total number of trucks required to travel between each origin-destination (O-D) pair to transport freight). Each private vehicle is driven by an individual road user who wants to minimize the personal travel cost/delay in the same time horizon. Thus, the two types of traffic compete for the limited road capacity. Since the individual road users do not cooperate, the impact of each private vehicle's travel behavior on the network traffic is not comparable to that of the truck company. So we model the problem as a Stackelberg game where the truck company is the leader and the other individual road users are the followers. We assume that private vehicles' travel behaviors follow dynamic user equilibrium given those of the trucks. This assumption may not hold for a single-day dynamic traffic assignment problem since it takes time to reach traffic equilibrium. However, we argue that the assumption is acceptable for a truck company that aims at making an optimal everyday transportation schedule in a long time horizon, especially when the travel demands are steady over time.

To understand the interaction of the two types of traffic in the Stackelberg game, we need to address the following questions: 1) How do the trucks' travel behaviors affect the private vehicles? 2) Considering the potential reactions of private vehicles towards its decision, how should the truck company schedule its freight transportation? In the remainder of this section, Section 3.1 and 3.2 answer the former question and Section 3.3 addresses the latter one.

### 3.1 Dynamic User Equilibrium with Inhomogeneous Traffic

We denote the time interval of interest by $\left[t_0, t_f\right] \subset \Re_+^1$. We assume that for all travelers (trucks and private vehicles) there is an identical desired time to arrive at their destinations, $T_A$. Note that this assumption can be relaxed by allowing different desired arrival time and we can just classify travelers with the same desired arrival time into the same group and model each group explicitly in the DUE which, although complicating the problem formulation, does not change the analysis and structure of the solution approach. To distinguish parameters and variables for trucks from those for private vehicles, we adopt superscript "*tr*" and "*pr*", respectively.

The unit cost to each traveler on a specific path is denoted by the effective unit path delay operator $\Psi_p\left(t, h^{pr}, h^{tr}\right)$, which is defined by the addition of travel cost and penalty cost:

$$\Psi_p\left(t, h^{pr}, h^{tr}\right) = D_p\left(t, h^{pr}, h^{tr}\right) + F\left[t + D_p\left(t, h^{pr}, h^{tr}\right) - T_A\right] \quad \forall p \in \mathrm{P},$$

where $D_p\left(t, h^{pr}, h^{tr}\right)$, $\forall p \in \mathrm{P}$ denotes the unit delay on path $p$ and $\mathrm{P}$ is the set of all paths employed by travelers. $t$ denotes the departure time and $h^{pr}, h^{tr}$ denote vectors of flows/departure rates (number of vehicles/trucks entering a path per unit time) of private vehicles and trucks, respectively. A more rigorous definition of $h$ can be found in Friesz (2010) and for simplicity is not stated here. $D_p\left(t, h^{pr}, h^{tr}\right)$ is an explicit function of the traveler's departure time but an implicit function of the other travelers' (including trucks' and private vehicles') travel behaviors. The quantification of $D_p\left(t, h^{pr}, h^{tr}\right)$ is discussed in detail in Section 3.2. $F\left[t + D_p\left(t, h^{pr}, h^{tr}\right) - T_A\right]$ denotes early/late arrival penalties and the following equations hold:

$$\begin{cases} \text{if} \quad t + D_p\left(t, h^{pr}, h^{tr}\right) > T_A, & F\left[t + D_p\left(t, h^{pr}, h^{tr}\right) - T_A\right] > 0 \\ \text{if} \quad t + D_p\left(t, h^{pr}, h^{tr}\right) < T_A, & F\left[t + D_p\left(t, h^{pr}, h^{tr}\right) - T_A\right] > 0 \\ \text{if} \quad t + D_p\left(t, h^{pr}, h^{tr}\right) = T_A, & F\left[t + D_p\left(t, h^{pr}, h^{tr}\right) - T_A\right] = 0 \end{cases}.$$

For simplicity, in this paper we assume a quadratic penalty function $F[x] = \alpha x^2$ where $\alpha$ is a constant.

To satisfy the travel demand for all private vehicles, the following flow conservation law must hold:

$$\sum_{p \in \mathrm{P}_{ij}} \int_{t_0}^{t_f} h_p^{pr}(t)\, dt = Q_{ij}^{pr} \quad \forall (i, j) \in W,$$

where $(i, j)$ is an O-D pair and $W$ is the set of all O-D pairs. $Q_{ij}^{pr} \in \Re_{++}^1$ is the fixed travel demand on O-D pair $(i, j)$ for private vehicles and $\mathrm{P}_{ij} \subset \mathrm{P}$ is a set of paths that connect O-D pair $(i, j)$.

The set of all feasible private vehicle flows can then be defined as

$$\Lambda_0 = \left\{ h_p^{pr} \geq 0 : \sum_{p \in P_{ij}} \int_{t_0}^{t_f} h_p^{pr}(t)\,dt = Q_{ij}^{pr} \quad \forall (i,j) \in W \right\} \subset \left( L_+^2 \left[ t_0, t_f \right] \right)^{|P|}.$$

**Theorem 1.** A vector of path flows $h^* \in \Lambda_0$ is a DUE if

$$h_p^*(t) > 0,\ p \in P_{ij} \Rightarrow \Psi_p\left[t, h^*(t)\right] = v_{ij}$$

where $v_{ij} = ess\inf\left[\Psi_p(t,h) : p \in P_{ij}\right], \quad \forall (i,j) \in W$.

**Proof.** See Friesz et al. (1993). □

Let's denote this equilibrium by $DUE\left(\Psi, \Lambda_0, t_0, t_f\right)$.

**Theorem 2.** $DUE\left(\Psi, \Lambda_0, t_0, t_f\right)$ is equivalent to the following differential variational inequality (DVI):

$$\left.\begin{array}{l} \text{find } h^* \in \Lambda \text{ such that} \\ \sum_{p \in P_{ij}} \int_{t_0}^{t_f} \Psi_p\left(t, h^*\right)\left(h_p - h_p^*\right) dt \geq 0 \\ \forall h \in \Lambda \end{array}\right\} DVI\left(\Psi, \Lambda, t_0, t_f\right)$$

where $\Lambda = \left\{ h \geq 0 : \frac{dy_{ij}}{dt} = \sum_{p \in P_{ij}} h_p(t),\ y_{ij}(t_0) = 0,\ y_{ij}\left(t_f\right) = Q_{ij} \quad \forall (i,j) \in W \right\}$ and $y_{ij}(t)$ stands for the volume of traffic that arrived at destination $j$ at time $t$.

**Proof.** See Friesz et al. (2011). □

In our model, given path flow of trucks, that of private vehicles will follow DUE. So the associated DVI formulation is

$$\left.\begin{array}{l} \text{find } h^{pr*} \in \Lambda \text{ such that} \\ \sum_{p \in P_{ij}} \int_{t_0}^{t_f} \Psi_p\left(t, h^{pr*}, h^{tr}\right)\left(h_p^{pr} - h_p^{pr*}\right) dt \geq 0 \\ \forall h^{pr} \in \Lambda \end{array}\right\} DVI^{pr}\left(\Psi, \Lambda, t_0, t_f\right) \tag{1}$$

where $\Lambda = \left\{ h^{pr} \geq 0 : \frac{dy_{ij}^{pr}}{dt} = \sum_{p \in P_{ij}} h_p^{pr}(t),\ y_{ij}^{pr}(t_0) = 0,\ y_{ij}^{pr}\left(t_f\right) = Q_{ij}^{pr} \quad \forall (i,j) \in W \right\}$ and $y_{ij}^{pr}(t)$ stands for the volume of private vehicles that arrived at destination $j$ at time $t$.

To solve (1), the evaluation of $\Psi_p\left(t, h^{pr}, h^{tr}\right)$ given $h^{pr}$ and $h^{tr}$ is necessary and it requires the quantification of $D_p\left(t, h^{pr}, h^{tr}\right)$ through dynamic network loading (DNL) which can be based on a link delay model (LDM).

### 3.2 Link Delay Model Based Dynamic Network Loading with Inhomogeneous Traffic

Dynamic network loading (DNL) refers to "the determination of arc-specific volumes, arc-specific exit rates and experienced path delay when departure rates are known for each path" (Friesz et al. 2011). With this process, DUE can be solved efficiently.

Additional notation is necessary. Each path can be represented by a sequence of connected arcs $p = \left\{a_1, a_2, \ldots, a_i, \ldots, a_{m(p)}\right\}$ where $m(p)$ denotes the total number of arcs contained in path $p$. Flows of trucks and private vehicles that travel along path $p$ exiting arc $a_i$ are denoted by $g_{a_i}^{p,tr}$ and $g_{a_i}^{p,pr}$, respectively. Volume of trucks and private vehicles traveling along path $p$ on arc $a_i$ are denoted by $x_{a_i}^{p,tr}$ and $x_{a_i}^{p,pr}$, respectively.

It is straightforward to derive the following arc dynamics:

$$\left.\begin{aligned}
&\frac{dx_{a_1}^{p,pr}(t)}{dt} = h_p^{pr}(t) - g_{a_1}^{p,pr}(t) \quad \forall p \in \mathrm{P} \\
&\frac{dx_{a_i}^{p,pr}(t)}{dt} = g_{a_{i-1}}^{p,pr}(t) - g_{a_i}^{p,pr}(t) \quad \forall p \in \mathrm{P}, i \in [2, m(p)] \\
&x_{a_i}^{p,pr}(t_0) = x_{a_i 0}^{p,pr} \quad \forall p \in \mathrm{P}, i \in [1, m(p)] \\
&\frac{dx_{a_1}^{p,tr}(t)}{dt} = h_p^{tr}(t) - g_{a_1}^{p,tr}(t) \quad \forall p \in \mathrm{P} \\
&\frac{dx_{a_i}^{p,tr}(t)}{dt} = g_{a_{i-1}}^{p,tr}(t) - g_{a_i}^{p,tr}(t) \quad \forall p \in \mathrm{P}, i \in [2, m(p)] \\
&x_{a_i}^{p,tr}(t_0) = x_{a_i 0}^{p,tr} \quad \forall p \in \mathrm{P}, i \in [1, m(p)]
\end{aligned}\right\} \tag{2}$$

where $x_{a_i 0}^{p,pr}$ and $x_{a_i 0}^{p,tr}$ are initial volume of private vehicles and trucks traveling along path $p$ on arc $a_i$.

Let's define

$$\delta_{a_i p} = \begin{cases} 1, & \text{if arc } a_i \text{ belongs to path p} \\ 0, & \text{otherwise} \end{cases}.$$

The total private vehicle volume on a specific arc $a_i$ at time *t* is

$$x_{a_i}^{pr}(t) = \sum_{p \in \mathrm{P}_{ij}} \delta_{a_i p} x_{a_i}^{p,pr}(t).$$

The total truck volume on a specific arc $a_i$ at time *t* is

$$x_{a_i}^{tr}(t) = \sum_{p \in \mathrm{P}_{ij}} \delta_{a_i p} x_{a_i}^{p,tr}(t).$$

The unit path delay equals the accumulated delays on arcs that compose the path. We denote the delay that a vehicle will experience when it enters arc $a_i$ at time $t$ by $D_{a_i}\left(x_{a_i}(t)\right)$, which is a function of the traffic volume on arc $a_i$ at time *t*. For simplicity, we assume a linear delay function, that is, $D_{a_i}(x) = A_{a_i} + B_{a_i} x$ where $A_{a_i}$ and $B_{a_i}$ are two positive constants. Considering the fact that the size of a truck may be different from that of a private vehicle thus may have a greater impact on congestion than a private vehicle, we assume that

$$x_{a_i}(t) = x_{a_i}^{pr}(t) + \beta x_{a_i}^{tr}(t)$$

where $\beta$ is a constant and is greater than or equal to 1.

If a vehicle enters path $p$ at time $t$, we denote its exit time from arc $a_i$ of path $p$ by $\tau_{a_i}^p(t)$ and we have

$$\begin{gathered}
\tau_{a_1}^p(t) = t + D_{a_1}\left[x_{a_1}(t)\right] \quad \forall p \in P \\
\tau_{a_i}^p(t) = \tau_{a_{i-1}}^p(t) + D_{a_i}\left[x_{a_i}\left(\tau_{a_{i-1}}^p(t)\right)\right] \quad \forall p \in P, i \in [2, m(p)]
\end{gathered}.$$

By applying the chain rule, we can derive the following flow propagation constraints based on the above exit time functions (see Friesz et al. (2001) for more details):

$$\left.\begin{aligned}
&g_{a_1}^{p,pr}\left(t+D_{a_1}\left[x_{a_1}(t)\right]\right)\left(1+D'_{a_1}\left[x_{a_1}(t)\right]\dot{x}_{a_1}\right)=h_p^{pr}(t) \quad \forall p\in\mathrm{P}\\
&g_{a_1}^{p,tr}\left(t+D_{a_1}\left[x_{a_1}(t)\right]\right)\left(1+D'_{a_1}\left[x_{a_1}(t)\right]\dot{x}_{a_1}\right)=h_p^{tr}(t) \quad \forall p\in\mathrm{P}\\
&g_{a_i}^{p,pr}\left(t+D_{a_i}\left[x_{a_i}(t)\right]\right)\left(1+D'_{a_i}\left[x_{a_i}(t)\right]\left(\dot{x}_{a_i}\right)\right)=g_{a_{i-1}}^{p,pr}(t) \quad \forall p\in\mathrm{P}, i\in\left[2,m(p)\right]\\
&g_{a_i}^{p,tr}\left(t+D_{a_i}\left[x_{a_i}(t)\right]\right)\left(1+D'_{a_i}\left[x_{a_i}(t)\right]\left(\dot{x}_{a_i}\right)\right)=g_{a_{i-1}}^{p,tr}(t) \quad \forall p\in\mathrm{P}, i\in\left[2,m(p)\right]
\end{aligned}\right\}$$

where as defined earlier $x_{a_i}(t)=\sum_{p\in\mathrm{P}_{ij}}\delta_{a_i p}\left(x_{a_i}^{p,pr}(t)+\beta x_{a_i}^{p,tr}(t)\right)$. In the above flow propagation constraints, the superscript "⇔" denotes differentiation with respect to the associated function argument and the overdot "·" refers to differentiation with respect to time.

The time shifts ($t+D_{a_i}[\cdot]$) in the above flow propagation constraints complicate the computation of the solution of DNL. We adopt the approximation proposed by Friesz et al. (2011) to simplify the network loading procedure. By defining $r_{a_i}^{p}(t)=\dfrac{dg_{a_i}^{p}(t)}{dt} \quad \forall p\in\mathrm{P},\, i\in\left[1,m(p)\right]$, the flow propagation constraints can be approximated by

$$\left.\begin{aligned}
&\frac{dr_{a_1}^{p,pr}(t)}{dt}=R_{a_1}^{p,pr}(x,g,r,h) \quad \forall p\in\mathrm{P}\\
&\frac{dr_{a_i}^{p,pr}(t)}{dt}=R_{a_i}^{p,pr}(x,g,r) \quad \forall p\in\mathrm{P}, i\in\left[2,m(p)\right]\\
&\frac{dr_{a_1}^{p,tr}(t)}{dt}=R_{a_1}^{p,tr}(x,g,r,h) \quad \forall p\in\mathrm{P}\\
&\frac{dr_{a_i}^{p,tr}(t)}{dt}=R_{a_i}^{p,tr}(x,g,r) \quad \forall p\in\mathrm{P}, i\in\left[2,m(p)\right]
\end{aligned}\right\} \qquad (3)$$

where

$$R_{a_1}^{p,pr}(x,g,r,h)=\frac{2h_p^{pr}(t)}{\left(D_{a_1}\left[x_{a_1}(t)\right]\right)^2\left(1+D'_{a_1}\left[x_{a_1}(t)\right]\dot{x}_{a_1}\right)}-\frac{2\left(g_{a_1}^{p,pr}(t)+r_{a_1}^{p,pr}(t)D_{a_1}\left[x_{a_1}(t)\right]\right)}{\left(D_{a_1}\left[x_{a_1}(t)\right]\right)^2}$$
$$\forall p\in\mathrm{P}$$

$$R_{a_i}^{p,pr}(x,g,r)=\frac{2g_{a_{i-1}}^{p,pr}(t)}{\left(D_{a_i}\left[x_{a_i}(t)\right]\right)^2\left(1+D'_{a_i}\left[x_{a_i}(t)\right]\dot{x}_{a_i}\right)}-\frac{2\left(g_{a_i}^{p,pr}(t)+r_{a_i}^{p,pr}(t)D_{a_i}\left[x_{a_i}(t)\right]\right)}{\left(D_{a_i}\left[x_{a_i}(t)\right]\right)^2}$$
$$\forall p\in\mathrm{P}, i\in\left[2,m(p)\right]$$

$$R_{a_1}^{p,tr}(x,g,r,h)=\frac{2h_p^{tr}(t)}{\left(D_{a_1}\left[x_{a_1}(t)\right]\right)^2\left(1+D'_{a_1}\left[x_{a_1}(t)\right]\dot{x}_{a_1}\right)}-\frac{2\left(g_{a_1}^{p,tr}(t)+r_{a_1}^{p,tr}(t)D_{a_1}\left[x_{a_1}(t)\right]\right)}{\left(D_{a_1}\left[x_{a_1}(t)\right]\right)^2}$$
$$\forall p\in\mathrm{P}$$

$$R_{a_i}^{p,tr}(x,g,r)=\frac{2g_{a_{i-1}}^{p,tr}(t)}{\left(D_{a_i}\left[x_{a_i}(t)\right]\right)^2\left(1+D'_{a_i}\left[x_{a_i}(t)\right]\dot{x}_{a_i}\right)}-\frac{2\left(g_{a_i}^{p,tr}(t)+r_{a_i}^{p,tr}(t)D_{a_i}\left[x_{a_i}(t)\right]\right)}{\left(D_{a_i}\left[x_{a_i}(t)\right]\right)^2}$$
$$\forall p\in\mathrm{P}, i\in\left[2,m(p)\right]$$

**Theorem 3**. If (2) and (3) satisfy the following regularity conditions:

(a). the arc exit time functions $\tau_{a_i}^{p}(t)$ for all *i* are strictly monotonic;

(b). the arc delay functions $D_{a_i}\left(x_{a_i}(t)\right)$ for all *i* are bounded and strictly positive;

(c). $D''_{a_i}\left(x_{a_i}(t)\right)$ exists and is continuous;

(d). $h_p^{tr}(t)$ and $h_p^{pr}(t)$ are continuous.

there exists a unique solution to (2) and (3).

**Proof.** Per Walter (1988), a unique solution to (2) and (3) exists if the right-hand side of each differential equation in (2) and (3) is continuously differentiable with respect to $x_{a_i}^{pr}$, $x_{a_i}^{tr}$, $g_{a_i}^{p,pr}$, $g_{a_i}^{p,tr}$, $r_{a_i}^{p,pr}$ and $r_{a_i}^{p,tr}$. It is obvious that the differential equations in (2) have this property. The remaining effort is on proving that those in (3) have the same property, that is, $R_{a_i}^{p,pr}$, $R_{a_i}^{p,tr}$ are differentiable with respect to $x_{a_i}^{pr}$, $x_{a_i}^{tr}$, $g_{a_i}^{p,pr}$, $g_{a_i}^{p,tr}$, $r_{a_i}^{p,pr}$ and $r_{a_i}^{p,tr}$. For simplicity, here we just show the closed-form expressions for those derivatives of $R_{a_1}^{p,pr}$ and $R_{a_1}^{p,tr}$, similar results can be established for $R_{a_i}^{p,pr}$ and $R_{a_i}^{p,tr}$ for $i\in\left[2,m(p)\right]$.

It is not difficult to derive the following closed-form expressions for derivatives of $R_{a_i}^{p,pr}$ and $R_{a_i}^{p,tr}$:

$$\frac{\partial R_{a_1}^{p,pr}}{\partial x_{a_1}^{pr}}=\frac{\partial R_{a_1}^{p,pr}}{\partial x_{a_1}^{tr}}=X_{a_1}^{pr}-\frac{2r_{a_1}^{p,pr}(t)D'_{a_1}\left[x_{a_1}(t)\right]}{\left(D_{a_1}\left[x_{a_1}(t)\right]\right)^2}+\frac{4\left(g_{a_1}^{p,pr}(t)+r_{a_1}^{p,pr}(t)D_{a_1}\left[x_{a_1}(t)\right]\right)D_{a_1}\left[x_{a_1}(t)\right]D'_{a_1}\left[x_{a_1}(t)\right]}{\left(D_{a_1}\left[x_{a_1}(t)\right]\right)^4}$$

$\forall p\in \mathrm{P}$

$$\frac{\partial R_{a_1}^{p,tr}}{\partial x_{a_1}^{pr}}=\frac{\partial R_{a_1}^{p,tr}}{\partial x_{a_1}^{tr}}=X_{a_1}^{tr}-\frac{2r_{a_1}^{p,tr}(t)D'_{a_1}\left[x_{a_1}(t)\right]}{\left(D_{a_1}\left[x_{a_1}(t)\right]\right)^2}+\frac{4\left(g_{a_1}^{p,tr}(t)+r_{a_1}^{p,tr}(t)D_{a_1}\left[x_{a_1}(t)\right]\right)D_{a_1}\left[x_{a_1}(t)\right]D'_{a_1}\left[x_{a_1}(t)\right]}{\left(D_{a_1}\left[x_{a_1}(t)\right]\right)^4}$$

$\forall p\in \mathrm{P}$

$$\frac{\partial R_{a_1}^{p,pr}}{\partial g_{a_1}^{p,pr}}=\frac{2h_p^{pr}(t)}{\left(D_{a_1}\left[x_{a_1}(t)\right]\right)^2}\frac{D'_{a_1}\left[x_{a_1}(t)\right]}{\left(1+D'_{a_1}\left[x_{a_1}(t)\right]\left(h_p^{pr}(t)-g_{a_1}^{p,pr}(t)\right)\right)^2}-\frac{2g_{a_1}^{p,pr}(t)}{\left(D_{a_1}\left[x_{a_1}(t)\right]\right)^2}\quad \forall p\in \mathrm{P}$$

$$\frac{\partial R_{a_1}^{p,tr}}{\partial g_{a_1}^{p,tr}}=\frac{2h_p^{tr}(t)}{\left(D_{a_1}\left[x_{a_1}(t)\right]\right)^2}\frac{D'_{a_1}\left[x_{a_1}(t)\right]}{\left(1+D'_{a_1}\left[x_{a_1}(t)\right]\left(h_p^{tr}(t)-g_{a_1}^{p,tr}(t)\right)\right)^2}-\frac{2g_{a_1}^{p,tr}(t)}{\left(D_{a_1}\left[x_{a_1}(t)\right]\right)^2}\quad \forall p\in \mathrm{P}$$

$$\frac{\partial R_{a_1}^{p,pr}}{\partial r_{a_1}^{p,pr}}=\frac{\partial R_{a_1}^{p,tr}}{\partial r_{a_1}^{p,tr}}-\frac{2}{D_{a_1}\left[x_{a_1}(t)\right]}\quad \forall p\in \mathrm{P}$$

where

$$X_{a_1}^{pr} = \frac{-2h_p^{pr}(t)\left(Y_{a_1}^{pr} + Z_{a_1}^{pr}\right)}{\left(D_{a_1}\left[x_{a_1}(t)\right]\right)^4\left(1 + D'_{a_1}\left[x_{a_1}(t)\right]\left(h_p^{pr}(t) - g_{a_1}^{p,pr}(t)\right)\right)^2} \quad \forall p \in \mathrm{P}$$

$$X_{a_1}^{tr} = \frac{-2h_p^{tr}(t)\left(Y_{a_1}^{tr} + Z_{a_1}^{tr}\right)}{\left(D_{a_1}\left[x_{a_1}(t)\right]\right)^4\left(1 + D'_{a_1}\left[x_{a_1}(t)\right]\left(ht_p^{pr}(t) - g_{a_1}^{p,tr}(t)\right)\right)^2} \quad \forall p \in \mathrm{P}$$

$$Y_{a_1}^{pr} = 2D_{a_1}\left[x_{a_1}(t)\right]D'_{a_1}\left[x_{a_1}(t)\right]\left(1 + D'_{a_1}\left[x_{a_1}(t)\right]\left(h_p^{pr}(t) - g_{a_1}^{p,pr}(t)\right)\right) \quad \forall p \in \mathrm{P}$$

$$Y_{a_1}^{tr} = 2D_{a_1}\left[x_{a_1}(t)\right]D'_{a_1}\left[x_{a_1}(t)\right]\left(1 + D'_{a_1}\left[x_{a_1}(t)\right]\left(h_p^{tr}(t) - g_{a_1}^{p,tr}(t)\right)\right) \quad \forall p \in \mathrm{P}$$

$$Z_{a_1}^{pr} = \left(D_{a_1}\left[x_{a_1}(t)\right]\right)^2\left(D''_{a_1}\left[x_{a_1}(t)\right]\left(h_p^{pr}(t) - g_{a_1}^{p,pr}(t)\right)\right) \quad \forall p \in \mathrm{P}$$

$$Z_{a_1}^{tr} = \left(D_{a_1}\left[x_{a_1}(t)\right]\right)^2\left(D''_{a_1}\left[x_{a_1}(t)\right]\left(h_p^{tr}(t) - g_{a_1}^{p,tr}(t)\right)\right) \quad \forall p \in \mathrm{P}$$

Thus the theorem is proven. □

So given $h^{pr}$ and $h^{tr}$ as constants, solving (2) and (3), also known as the differential algebraic equation (DAE) system, we can evaluate $x_{a_i}(t)$, $D_{a_i}\left[x_{a_i}(t)\right]$ and $\tau_{a_i}^p(t)$ for all $p \in \mathrm{P}, i \in \left[2, m(p)\right]$, and then quantify unit path delay using the following function

$$D_p(t) = \sum_{i=1}^{m(p)}\left[\tau_{a_i}^p(t) - \tau_{a_{i-1}}^p(t)\right] = \tau_{a_{m(p)}}^p(t) - t .$$

It was first proven by Friesz et al. (1993) that a closed-form formulation of DUE can preserve first-in-first-out (FIFO) rule under mild assumptions. In the remaining of this section, we prove a similar result for DUE with inhomogeneous traffic.

First we need to introduce the following lemma.

**Lemma 1**. If function $f : \Re \to \Re$ is invertible and differentiable with a derivative $f'$:

$$\left[f^{-1}\right]'(z) = 1 / f'\left[f^{-1}(z)\right]$$

**Proof.** See Friesz et al. (1993).

Now we can demonstrate the following theorem.

**Theorem 4**. For an arc delay function $D_{a_i}\left(x_{a_i}(t)\right) = A_{a_i} + B_{a_i}\left(x_{a_i}^{pr}(t) + \beta x_{a_i}^{tr}(t)\right)$, the resulting arc exit time function $\tau_{a_i}$ is strictly increasing thus invertible. Consequently, the FIFO rule is satisfied.

**Proof.** We partition the time into appropriate intervals and assume without loss of generality that $x_{a_i}^{pr}(0) = x_{a_i}^{tr}(0) = 0$. We in addition assume that the first vehicle enters arc $a_i$ at time 0 and it does not matter whether it is a private vehicle or a truck. We denote the time that the first vehicle exits arc $a_i$ by $t_1$ and we have by definition

$$t_1 = D_{a_i}(0) = A_{a_i}$$

and

$$x_{a_i}(t) = x_{a_i}^{pr}(t) + \beta x_{a_i}^{tr}(t) = \int_0^t u_{a_i}^{pr}(s)ds + \beta\int_0^t u_{a_i}^{tr}(s)ds, \quad \forall t \in [0, t_1]$$

where $u_{a_i}^{pr}(t)=\sum_{p\in P} h_p^{pr}(t)\delta_{a_i p}$ and $u_{a_i}^{tr}(t)=\sum_{p\in P} h_p^{tr}(t)\delta_{a_i p}$ denote the flow of private vehicles and trucks entering arc $a_i$ , respectively.

Hence for any $t\in[0,t_1]$, we denote the exit time function by $\tau_{1,a_i}(t)$, and we have

$$\tau_{1,a_i}(t)=t+D_{a_i}(t)=t+A_{a_i}+B_{a_i}\left[\int_0^t u_{a_i}^{pr}(s)ds+\beta\int_0^t u_{a_i}^{tr}(s)ds\right].$$

Note that $t_1=\tau_{1,a_i}(0)$ and $\tau_{1,a_i}(t)$ is differentiable with

$$\tau'_{1,a_i}(t)=1+B_{a_i}\left[u_{a_i}^{pr}(t)+\beta u_{a_i}^{tr}(t)\right] \tag{i}$$

which is strictly positive since $B_{a_i}$ is positive and $u_{a_i}^{pr}(t)$ and $u_{a_i}^{tr}(t)$ are nonnegative. Hence, $\tau_{1,a_i}(t)$ is increasing on $[0,t_1]$ and a well defined inverse function $\tau_{1,a_i}^{-1}(t)$ exists on $[\tau_{a_1}(0),\tau_{a_1}(t_1)]$.

Let $t_2=\tau_{a_1}(t_1)$, $[\tau_{a_1}(0),\tau_{a_1}(t_1)]\equiv[t_1,t_2]$.

Note that the commuters traveling on arc $a_i$ at time $t\in[t_1,t_2]$ are those who entered the arc during the interval $[\tau_{1,a_i}^{-1}(t),t]$. Hence, if we denote the exit time for commuters entering arc $a_i$ during the interval $[t_1,t_2]$ by $\tau_{2,a_i}(t)$, then for all $t\in[t_1,t_2]$,

$$\tau_{2,a_i}(t)=t+A_{a_i}+B_{a_i}\left[\int_{\tau_{1,a_i}^{-1}(t)}^t u_{a_i}^{pr}(s)ds+\beta\int_{\tau_{1,a_i}^{-1}(t)}^t u_{a_i}^{tr}(s)ds\right]$$

and $\tau_{2,a_i}(t_1)=\tau_{1,a_i}(t_1)=t_2$ .

The derivative of $\tau_{2,a_i}(t)$ is given by

$$\begin{aligned}
\tau'_{2,a_i}(t)&=1+B_{a_i}\left[u_{a_i}^{pr}(t)-u_{a_i}^{pr}\left[\tau_{1,a_i}^{-1}(t)\right]\left[\tau_{1,a_i}^{-1}\right]'(t)+\beta u_{a_i}^{tr}(t)-\beta u_{a_i}^{tr}\left[\tau_{1,a_i}^{-1}(t)\right]\left[\tau_{1,a_i}^{-1}\right]'(t)\right]\\
&=1+B_{a_i}u_{a_i}^{pr}(t)-\frac{B_{a_i}u_{a_i}^{pr}\left[\tau_{1,a_i}^{-1}(t)\right]}{\tau'_{1,a_i}\left[\tau_{1,a_i}^{-1}(t)\right]}+B_{a_i}\beta u_{a_i}^{tr}(t)-\frac{B_{a_i}\beta u_{a_i}^{tr}\left[\tau_{1,a_i}^{-1}(t)\right]}{\tau'_{1,a_i}\left[\tau_{1,a_i}^{-1}(t)\right]}\\
&=1+B_{a_i}u_{a_i}^{pr}(t)-\frac{B_{a_i}u_{a_i}^{pr}\left[\tau_{1,a_i}^{-1}(t)\right]}{1+B_{a_i}\left(u_{a_i}^{pr}\left[\tau_{1,a_i}^{-1}(t)\right]+\beta u_{a_i}^{tr}\left[\tau_{1,a_i}^{-1}(t)\right]\right)}+B_{a_i}\beta u_{a_i}^{tr}(t)-\frac{B_{a_i}\beta u_{a_i}^{tr}\left[\tau_{1,a_i}^{-1}(t)\right]}{1+B_{a_i}\left(u_{a_i}^{pr}\left[\tau_{1,a_i}^{-1}(t)\right]+\beta u_{a_i}^{tr}\left[\tau_{1,a_i}^{-1}(t)\right]\right)}\\
&=B_{a_i}\left(u_{a_i}^{pr}(t)+\beta u_{a_i}^{tr}(t)\right)-\frac{1}{1+B_{a_i}\left(u_{a_i}^{pr}\left[\tau_{1,a_i}^{-1}(t)\right]+\beta u_{a_i}^{tr}\left[\tau_{1,a_i}^{-1}(t)\right]\right)}\\
&>B_{a_i}\left(u_{a_i}^{pr}(t)+\beta u_{a_i}^{tr}(t)\right)\geq 0
\end{aligned}$$

The second and third equalities are derived using Lemma 1 and equation (i), respectively.

Thus, we can again conclude that $\tau_{2,a_i}(t)$ is increasing on $[t_1,t_2]$ and a well defined inverse function $\tau_{2,a_i}^{-1}(t)$ exists on $[\tau_{a_1}(t_1),\tau_{a_1}(t_2)]$.

We proceed by induction. For *n = 2*, we have already shown that

$$\tau_{2,a_i}(t)=t+A_{a_i}+B_{a_i}\left[\int_{\tau_{1,a_i}^{-1}(t)}^t u_{a_i}^{pr}(s)ds+\beta\int_{\tau_{1,a_i}^{-1}(t)}^t u_{a_i}^{tr}(s)ds\right]$$

$$\tau_{2,a_i}(t_1)=\tau_{1,a_i}(t_1)=t_2$$

$$\tau'_{2,a_i}(t)>B_{a_i}\left(u_{a_i}^{pr}(t)+\beta u_{a_i}^{tr}(t)\right)\geq 0$$

Then we choose any $k>2$ . Suppose that for $n=k$ , the exit time function $\tau_{k,a_i}(t)$ is invertible and the following conditions hold:

$$\tau_{k,a_i}(t) = t + A_{a_i} + B_{a_i}\left[\int_{\tau^{-1}_{k-1,a_i}(t)}^{t} u^{pr}_{a_i}(s)ds + \beta\int_{\tau^{-1}_{k-1,a_i}(t)}^{t} u^{tr}_{a_i}(s)ds\right] \quad \text{(ii)}$$

$$\tau_{k,a_i}(t_{k-1}) = \tau_{k-1,a_i}(t_{k-1}) = t_k \quad \text{(iii)}$$

$$\tau'_{k,a_i}(t) > B_{a_i}\left(u^{pr}_{a_i}(t) + \beta u^{tr}_{a_i}(t)\right) \text{ for all } t \in [t_{k-1}, t_k] \quad \text{(iv)}$$

We wish to show that for *n = k + 1*, if denote $t_{k+1} = \tau_{k,a_i}(t_k)$, then the exit time function $\tau_{k+1,a_i}(t_k)$ satisfies the following conditions:

$$\tau_{k+1,a_i}(t) = t + A_{a_i} + B_{a_i}\left[\int_{\tau^{-1}_{k,a_i}(t)}^{t} u^{pr}_{a_i}(s)ds + \beta\int_{\tau^{-1}_{k,a_i}(t)}^{t} u^{tr}_{a_i}(s)ds\right] \quad \text{(v)}$$

$$\tau_{k+1,a_i}(t_k) = \tau_{k,a_i}(t_k) = t_{k+1} \quad \text{(vi)}$$

$$\tau'_{k+1,a_i}(t) > B_{a_i}\left(u^{pr}_{a_i}(t) + \beta u^{tr}_{a_i}(t)\right) \text{ for all } t \in [t_k, t_{k+1}] \quad \text{(vii)}$$

By definition, it is not difficult to show that equation (v) is satisfied. So we only need to derive equations (vi) and (viii) in the rest of this proof.

First, similar to the derivation of $\tau'_{2,a_i}(t)$ we can show equation (vii) by the following:

$$\begin{aligned}
\tau'_{k+1,a_i}(t) &= 1 + B_{a_i}\left[u^{pr}_{a_i}(t) - u^{pr}_{a_i}\left[\tau^{-1}_{k,a_i}(t)\right]\left[\tau^{-1}_{k,a_i}\right]'(t) + \beta u^{tr}_{a_i}(t) - \beta u^{tr}_{a_i}\left[\tau^{-1}_{k,a_i}(t)\right]\left[\tau^{-1}_{k,a_i}\right]'(t)\right] \\
&= 1 + B_{a_i}u^{pr}_{a_i}(t) - \frac{B_{a_i}u^{pr}_{a_i}\left[\tau^{-1}_{k,a_i}(t)\right]}{\tau'_{k,a_i}\left[\tau^{-1}_{k,a_i}(t)\right]} + B_{a_i}\beta u^{tr}_{a_i}(t) - \frac{B_{a_i}\beta u^{tr}_{a_i}\left[\tau^{-1}_{k,a_i}(t)\right]}{\tau'_{k,a_i}\left[\tau^{-1}_{k,a_i}(t)\right]} \\
&= B_{a_i}\left(u^{pr}_{a_i}(t) + \beta u^{tr}_{a_i}(t)\right) - \frac{1}{1 + B_{a_i}\left(u^{pr}_{a_i}\left[\tau^{-1}_{k,a_i}(t)\right] + \beta u^{tr}_{a_i}\left[\tau^{-1}_{k,a_i}(t)\right]\right)} \\
&> B_{a_i}\left(u^{pr}_{a_i}(t) + \beta u^{tr}_{a_i}(t)\right) \geq 0
\end{aligned}$$

Then by equation (iii) and the assumed invertibility of $\tau_{k,a_i}(t)$ and $\tau_{k-1,a_i}(t)$, we have

$$t_{k-1} = \tau^{-1}_{k,a_i}(t_k) \text{ and } t_{k-1} = \tau^{-1}_{k-1,a_i}(t_{k-1}).$$

This together with equation (v) yields

$$\begin{aligned}
\tau_{k+1,a_i}(t_k) &= t_k + A_{a_i} + B_{a_i}\left[\int_{\tau^{-1}_{k,a_i}(t_k)}^{t_k} u^{pr}_{a_i}(s)ds + \beta\int_{\tau^{-1}_{k,a_i}(t_k)}^{t_k} u^{tr}_{a_i}(s)ds\right] \\
&= t_k + A_{a_i} + B_{a_i}\left[\int_{t_{k-1}}^{t_k} u^{pr}_{a_i}(s)ds + \beta\int_{t_{k-1}}^{t_k} u^{tr}_{a_i}(s)ds\right] \\
&= t_k + A_{a_i} + B_{a_i}\left[\int_{\tau^{-1}_{k-1,a_i}(t_{k-1})}^{t_k} u^{pr}_{a_i}(s)ds + \beta\int_{\tau^{-1}_{k-1,a_i}(t_{k-1})}^{t_k} u^{tr}_{a_i}(s)ds\right] \\
&= \tau_{k,a_i}(t_k)
\end{aligned}$$

and thus equation (vi) holds.

Consequently, the exit time function $\tau_{a_i}(t)$ is everywhere continuous and increasing on $\Re_+$ which indicates that $\tau_{a_i}(t_1) > \tau_{a_i}(t_2)$ if $t_1 > t_2$. Thus the FIFO rule holds. □

### 3.3 Formulation of the Stackelberg Game

In the proposed Stackelberg game, the leader (the truck company) aims at minimizing its total cost while satisfying the travel demand over the time frame of interest, which can be represented by (4) and (5), respectively:

$$\min \sum_{p \in \mathrm{P}_{ij}} \int_{t_0}^{t_f} \Psi_p\left(t, h^{pr*}, h^{tr}\right) h^{tr}_p dt \quad (4)$$

$$\sum_{p \in \mathrm{P}_{ij}} \int_{t_0}^{t_f} h_p^{tr}(t)\,dt = Q_{ij}^{tr},\; h_p^{tr} \ge 0 \quad \forall (i,j) \in W \tag{5}$$

where

$$\Psi_p\left(t, h^{pr}, h^{tr}\right) = D_p\left(t, h^{pr}, h^{tr}\right) + F\left[t + D_p\left(t, h^{pr}, h^{tr}\right) - T_A\right] \quad \forall p \in \mathrm{P}.$$

Then the problem can be represented by {(1), (2), (3), (4), (5)} which is a dynamic mathematical program with equilibrium constraints (MPEC).

In (4), $h^{pr*}$ stands for the vector of equilibrium flows of private vehicles and $h^{tr}$ is the vector of truck flows. Equation (4) minimizes the total effective delay experienced by all the trucks that travel on the network in time interval $\left[t_0, t_f\right]$. Equation (5) is the demand consumption constraint for trucks. (4) and (5) compose the upper level of the dynamic MPEC. (1), (2) and (3) which characterize the equilibrium flows of private vehicles compose the lower level of the dynamic MPEC.

Existence of solutions to this dynamic MPEC is a crucial issue. It is resolved if we can show that for any set of truck flows $h_p^{tr}(t),\ \forall p$ that are feasible to (5), (1) has a solution, that is, $h^{pr*}$. If $\Psi_p\left(t, h^{pr}, h^{tr}\right)$ for all $p$ is continuous and the feasible set $\Lambda$ in (1) is compact, the fixed point theory of multi-valued mappings in topological vector spaces discussed by Browder (1968) can be applied to show that there exists a solution to (1). So, existence results are general if the DUE is based on the formulation (1). However, a rigorous proof of existence when regularity conditions are imposed entirely on the path delay operators without the assumption that the path flows are a priori bounded from above is still challenging and remains to be addressed by future studies.

## 4. ALGORITHM DESIGN

Although a rigorous proof is not available, we still have a general property of existence of solutions and we can develop or apply an algorithm to find the solution. Since the dynamic MPEC is nonconvex and the DNL process in solving the lower level DUE requires complex nested calculations, we reformulate the problem as a mathematical program with complementarity constraints (MPCC), which is a single-level problem computable by standard optimization techniques. In particular, the DUE is formulated as the following complementarity problem:

$$\begin{aligned} &\left(\Psi_p\left(t, h_p^{pr}, h_p^{tr}\right) - \mu_{ij}\right) \perp h_p^{pr} \quad \forall p \in ij \\ &\Psi_p\left(t, h_p^{pr}, h_p^{tr}\right) - \mu_{ij} \ge 0 \qquad \forall p \in ij \\ &h^{pr} \ge 0 \end{aligned} \tag{6}$$

where $h^{pr} \in \Lambda$, and $\mu := \left(\mu_{ij} : ij \in W\right)$. With complementarity constraints (6) substituted for DUE (1), the MPCC reformulation of the proposed urban freight model can be represented by {(2), (3), (4), (5), (6)}.

Since the complementarity constraints might not satisfy certain constraint qualifications that are necessary to guarantee convergence of the solution (Rodrigues and Monteiro, 2006), to solve the MPCC, we penalize the complementarity constraints and obtain the augmented objective function as:

$$\begin{aligned} &U\left(h^{pr}, h^{tr}, \mu, M\right) \\ &:= \sum\nolimits_{p \in \mathrm{P}_{ij}} \int_{t_0}^{t_f} \Psi_p\left(t, h^{pr}, h^{tr}\right) h_p^{tr}\, dt \\ &\quad + M \sum\nolimits_{p \in \mathrm{P}_{ij}} \int_{t_0}^{t_f} \left[\left(\Psi_p\left(t, h_p^{pr}, h_p^{tr}\right) - \mu_w\right) h_p^{pr}\right]^2 dt \\ &\quad + M \sum\nolimits_{p \in \mathrm{P}_{ij}} \int_{t_0}^{t_f} \left[\max\left\{\mu_{ij} - \Psi_p\left(t, h_p^{pr}, h_p^{tr}\right), 0\right\}\right]^2 dt \end{aligned} \tag{7}$$

where $M$ is a large number.

By substituting (7) for (4), we have finalized the urban freight transportation planning problem as a single-level nonlinear program {(2), (3), (5), (6), (7)}. We design the following projected gradient algorithm to solve this nonlinear program.

**Projected gradient algorithm:**

**Step 0**. Initialization. Identify an initial feasible solution $u^k := \left(h^{pr,k}, h^{tr,k}, \mu_{ij}^k\right)^T$ and set *k* = 0;

**Step 1**. Solve the optimal control subproblem:

$$\begin{aligned}
&\min_v \int_{t_0}^{t_f} \frac{1}{2}\left\| h^{pr,k}(t) - \alpha F_1\left(t, h^{pr,k}(t), h^{tr,k}(t), \mu^k\right) - h^{pr}(t) \right\|_2^2 dt \\
&+\int_{t_0}^{t_f} \frac{1}{2}\left\| h^{tr,k}(t) - \alpha F_2\left(t, h^{pr,k}(t), h^{tr,k}(t), \mu^k\right) - h^{tr}(t) \right\|_2^2 dt \\
&+\int_{t_0}^{t_f} \frac{1}{2}\left\| \mu^k - \alpha F_3\left(t, h^{pr,k}(t), h^{tr,k}(t), \mu^k\right) - \mu \right\|_2^2 dt \\
&s.t. \\
&\frac{dy_{ij}^{pr}(t)}{dt} = \sum_{p \in P_{ij}} h_p^{pr}(t) \\
&\frac{dy_{ij}^{tr}(t)}{dt} = \sum_{p \in P_{ij}} h_p^{tr}(t) \\
&y_{ij}^{pr}(t_0) = 0, y_{ij}^{pr}(t_f) = Q_{ij}^{pr} \\
&y_{ij}^{tr}(t_0) = 0, y_{ij}^{tr}(t_f) = Q_{ij}^{tr}
\end{aligned}$$

Denote the solution by $u^{k+1} := \left(h^{pr,k+1}, h^{tr,k+1}, \mu^{k+1}\right)^T$.

**Step 2**. Stop if $\left\| u^{k+1} - u^k \right\| \le \varepsilon_1$, where $\varepsilon_1 \in \Re_{++}^1$ is a preset scalar. Otherwise, set *M* = *CM*, go to Step 1.

In the pseudo-code, $F_1 := \partial U / \partial h^{pr}$, $F_2 := \partial U / \partial h^{tr}$, $F_3 := \partial U / \partial \mu$ and C are constants that are great than 1. Note that since the problem is nonconvex in general, by the proposed solution approach, only local optimal solutions can be guaranteed.

**5. NUMERICAL TEST**

The algorithm is tested on Nguyen-Dupuis network which consists of 19 nodes and 25 links as shown by Figure 1. We set $t_0 = 100, t_f = 175$ and define a desired arrival time $T_A = 160$ for the trucks as well as private vehicles. This planning horizon is then discretized into 300 time intervals and each of which is 0.25 time units long. Each time unit represents 1 minute in the real world.

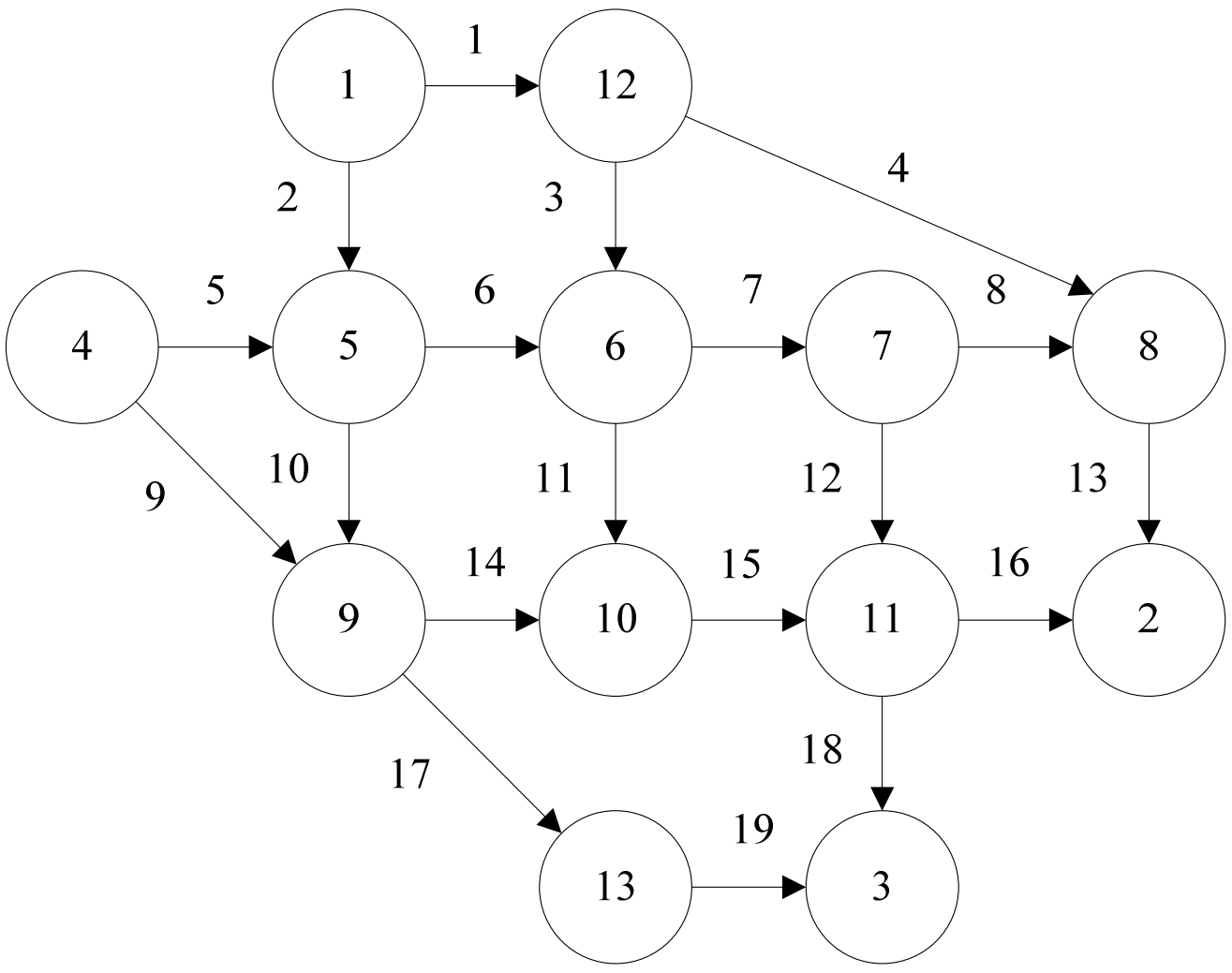

Figure 1 Nguyen-Dupuis network

Specifically, we focus on the set of O-D pairs $W = \{(1,2),(1,3),(4,2),(4,3)\}$ with each O-D pair having a fixed travel demand. Four test scenarios are considered. They are:

$$\begin{cases} \text{Scenario 1: } Q_{ij}^{pr} = 15000,\ Q_{ij}^{tr} = 5000 \\ \text{Scenario 2: } Q_{ij}^{pr} = 15000,\ Q_{ij}^{tr} = 2500 \\ \text{Scenario 3: } Q_{ij}^{pr} = 15000,\ Q_{ij}^{tr} = 500 \\ \text{Scenario 4: } Q_{ij}^{pr} = 15000,\ Q_{ij}^{tr} = 0 \end{cases}, \quad \forall (i,j) \in W .$$

We assume that in $D_{a_i}(x) = A_{a_i} + B_{a_i} x$, except that $A_{a_i} = 1.5$ for *i* = 1, 7, 13, 15 and 19, and $A_{a_i} = 2.5$ for *i* = 4, for remaining arcs $A_{a_i} = 2$. $B_{a_i} = 6.67 \times 10^{-4}$ for all arcs. We also assume that the trucks and private vehicles have the same length, that is, in $x_{a_i}(t) = x_{a_i}^{pr}(t) + \beta x_{a_i}^{tr}(t)$, $\beta = 1$. We set $\alpha = 0.5$ in early/late arrival penalties $F[x] = \alpha x^2$.

The solution approach is coded in MATLAB 7 and GAMS and solved on Penn State Lion-X system with the following attributes: Intel Xeon E5450 Quad-Core 3.0 GHz and 64GB RAM.

In the rest of this section, Sections 5.1, 5.2, and 5.3 address the following questions respectively.

- Will the interaction of freight and personal transportation influence the traffic flow pattern?
- To the truck company, what is the value of incorporating the interaction of freight and personal transportation?
- What managerial insight does the model bring to a MPO?

### 5.1 Traffic Flow Pattern

The "optimal" truck flows and equilibrium private vehicle flows are solved by the proposed algorithm. Since there are in total 25 paths that connect the 4 O-D pairs, without loss of generality we only analyze the time-varying traffic flows on Path 10 (Arcs 9, 14, 15 and 18 that connect O-D pair (4, 3)) in each scenario as illustrated in Figure 2.

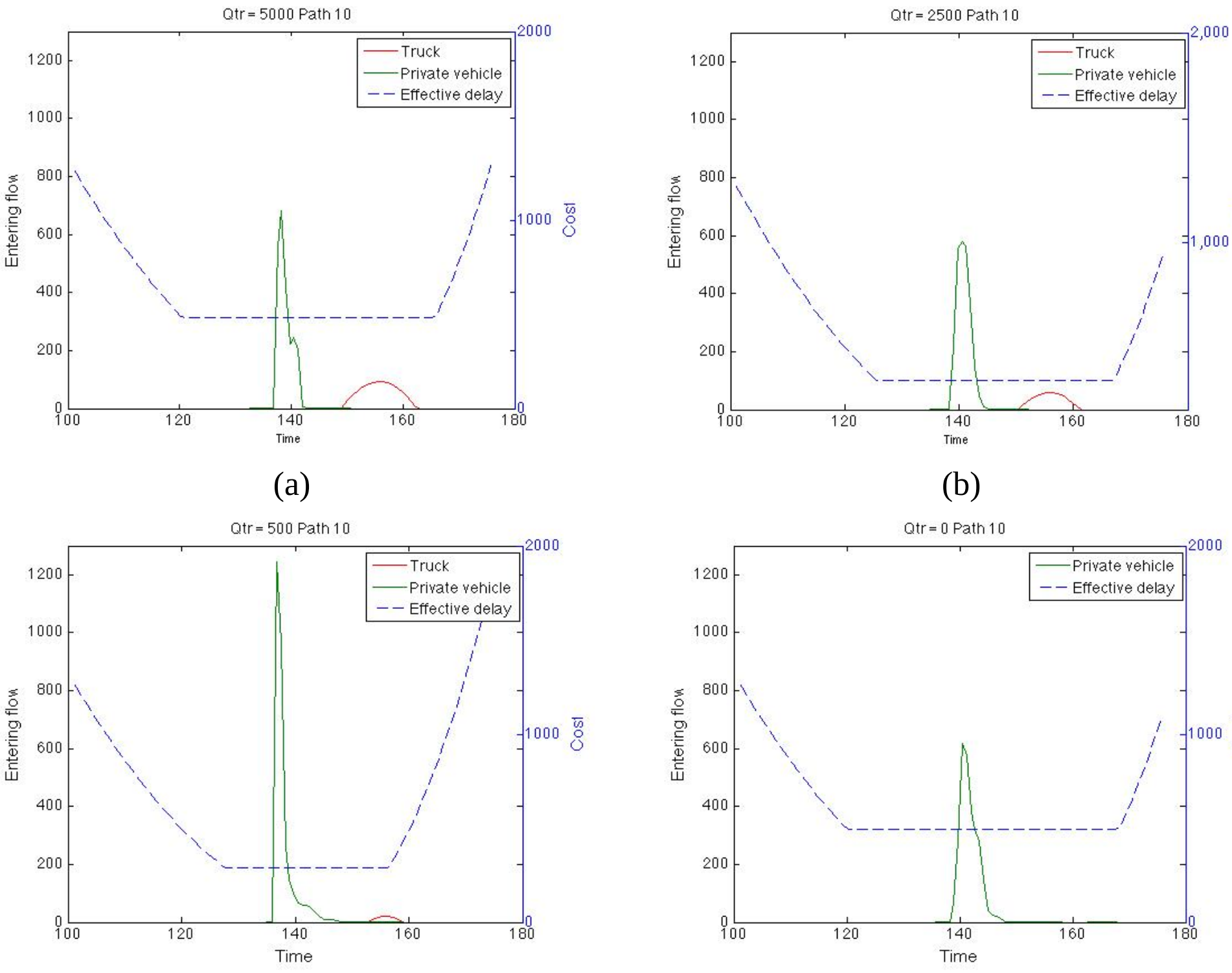

(a) (b)

(c) (d)

Figure 2 Dynamic traffic flows on Path 10 in 4 scenarios

In Figure 2, (a), (b), (c) and (d) show the dynamic traffic flows and the effective delay on Path 1 in Scenario 1, 2, 3 and 4, respectively. We can see that in all of these 4 scenarios, traffic flows exist only when the effective delay is at its minimum and the peak of the private vehicle flows occurs at around time 140 which can be considered as the appropriate departure time for private vehicles in order to arrive at the destination at the desired time $T_a$. From Figure 2, we can clearly identify the interaction of freight and personal transportation. As an example, in Scenario 3 (Figure 2(c)) when there are in total 500 trucks traveling on the network, the private vehicle flow significantly differs from that in Scenario 4 (Figure 2(d)) when there is no truck at all. What's more, again in Scenario 3, the peak of private vehicle flow occurs earlier that in Scenario 4. It is not intuitive why the traffic flow pattern in Scenario 3 significantly differs from those in Scenario 1 and 2. However, since the truck company optimizes its transportation cost over the entire network and Figure 2 just illustrates traffic flow pattern on one of the 25 paths of the network, such a result is not surprising.

**5.2 Comparative Study**

By Figure 2, we have demonstrated that there exist interactions between truck and private vehicles even when the number of trucks is relatively small. We would like to further estimate the importance of considering such an interaction in truck scheduling by considering two cases. In Case 1, the truck company ignores the interaction and it simply assumes that the private vehicles will follow DUE without taking into account the truck flows. It estimates the effective delay operator under the regular DUE and uses it as the cost coefficient to optimize truck flows. Let's mark the optimal solution in this case as $\left(\hat{h}_p^{tr}\right)$. However, as we know, the truck flow will influence the private vehicle flow so that given $\left(\hat{h}_p^{tr}\right)$, the real DUE flow for private vehicles is different from its estimation and thus the real effective delay operator is not the same as estimated. Thus in Case 1, the real total cost should be calculated using the real effective delay operator. Let's mark the real total cost as $z^0$. In Case 2, the truck company considers the interaction while optimizing truck flow and the optimal solution will be exactly the solution to the MPCC {(2), (3), (5), (6), (7)}. In Case 2, the minimal total cost is also the real total cost and let's mark it as $z^1$. We then compare $z^0$ and $z^1$ in Scenario 1, 2 and 3 in Table 1 which shows that in all 3 scenarios considering the interaction can help reduce the total cost to the truck company significantly, especially when freight transportation accounts for a considerable portion of total traffic. The computation time is summarized in Table 2 and it is not surprising that in Case 2 the computation time in each scenario is much greater than its counterpart in Case 1 since in Case 1 the original problem is just single-level optimization thus does not involve iterative updates of DUE for private vehicles. However, considering the significant reduction in cost, it is still worth incorporating the interaction of freight and personal transportation for the truck company while making truck schedules.

Table 1 Comparison of total cost

| | Scenario 1 | Scenario 2 | Scenario 3 |
|---|---|---|---|
| $z^0$ | $1.42\times10^7$ | $4.52\times10^6$ | $7.99\times10^5$ |
| $z^1$ | $2.58\times10^7$ | $7.30\times10^6$ | $9.09\times10^5$ |
| Reduction | 44.90% | 38.16% | 12.12% |

Table 2 Comparison of computation time (in second)

| | Scenario 1 | Scenario 2 | Scenario 3 |
|---|---|---|---|
| Case 1 | 2180.46 | 1760.42 | 1851.65 |
| Case 2 | 27408.29 | 29994.63 | 63580.18 |

**5.3 Braess-Like Paradox in Dynamic Stackelberg Game**

So far in this paper, we have demonstrated that the interaction of freight and personal transportation should be considered when a truck company schedules its freight transportation in order to save cost. In this subsection, we want to address the concern of a MPO whose objective is to minimize total social cost, the cost to the all owners of private vehicles and the truck company, while ensuring that the travel demand is satisfied. A possible approach is to make a policy to restrict trucks from entering a portion of the network. However, it is not necessarily reducing the total cost since the transportation cost to the truck company may soar. So if it works, we actually come across a Braess-like Paradox: reducing capacity of a network for partial road users who selfishly select their routes can increase overall performance (see Akamatsu and Heydecker (2003) and Lin and Lo (2009) for more details about the dynamic extension of Braess Paradox). Note that since the MPO knows exactly that the truck company and the owners of private vehicles play a Stackelberg game, the problem is actually a tri-level optimization problem. To simplify the formulation and solution of the problem, we explore the existence of such Braess-like paradoxes by a trial-and-error approach. We arbitrarily choose one arc at a time from the network to be blocked for trucks. Then we conduct numerical test to check whether blocking that arc for trucks can reduce the total cost: if yes, then we find the Braess-like paradox and we can stop; if no, we choose another arc and repeat the numerical test and the judgment of result. After several iterations, we find that blocking all the trucks from entering Arc 12 during the time interval of our interest can reduce the total cost (see Figure 3). More details of this policy are summarized in Table 3.

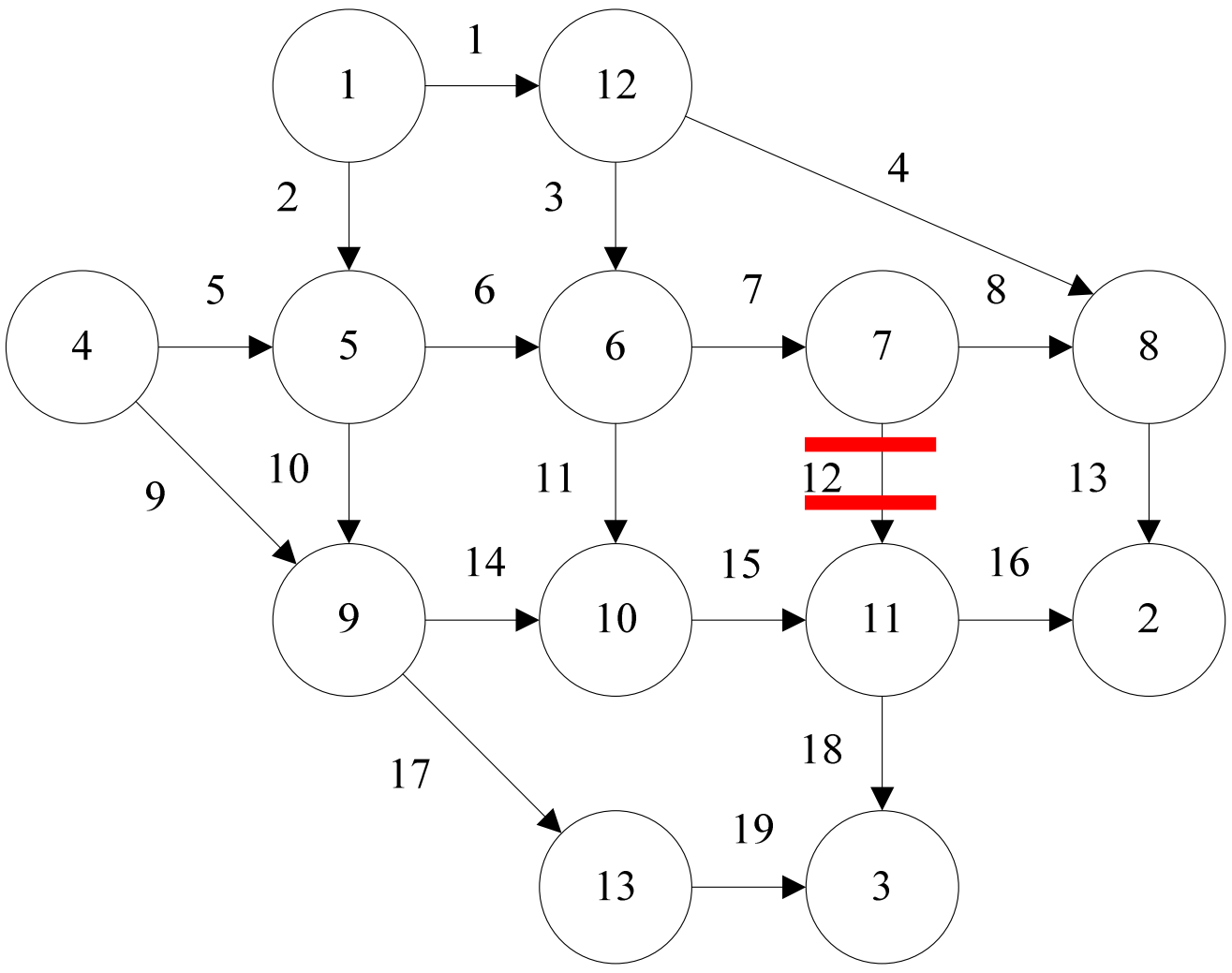

Figure 3 Nguyen-Dupuis network with Arc 12 blocked for all trucks

Table 3 Braess-like Paradox

| | Total Cost | Truck Company's Cost | Computation time (in second) |
|---|---|---|---|
| All arcs are accessible | $3.16\times10^7$ | $4.52\times10^6$ | 29994.63 |
| Arc 12 blocked for trucks | $3.11\times10^7$ | $5.02\times10^6$ | 23082.89 |
| Improvement | 1.6% | -11.1% | 23.0% |

From Table 3, we know that the MPO can block trucks from entering Arc 12 to reduce total cost by 1.6%. However, it may increase the truck company's cost by 11.1%. So before the MPO implement such a policy, it should consider the possible obstruction from the truck company and try to balance the increase in social welfare and the loss in the truck company's profit. As an example resolution, the MPO may compensate the truck company for restricting its travel right.

## 6. CONCLUSION

In this paper we model urban freight transportation planning by a dynamic Stackelberg game and formulate the problem as a dynamic MPEC. In particular, we assume that there is a truck company trying to minimize its freight transportation cost which is dependent on its transportation plan and the background traffic such as personal transportation. The model explicitly characterizes the interaction of freight and personal transportation by its lower level problem, that is, DUE with inhomogeneous traffic, and optimizes the truck schedule in the upper level which is an SO problem. To obtain local optimal solutions and achieve efficient computation, we reformulate the MPEC as a MPCC and design a projected gradient algorithm to solve it. Numerical results show that the interaction between different road users - specifically, a trucking company and individuals, exists and is nonnegligible even when the amount of trucks compared to that of private vehicles is small, and demonstrate that significant cost reduction can be achieved if the truck company schedules freight transportation considering this interaction, which supports our concern that in an urban network personal transportation should not be ignored while scheduling freight transportation. Moreover, with extensive numerical tests we find a Braess-like paradox in this dynamic bi-level transportation planning problem which implicates that a MPO may increase social welfare by restricting trucks from entering specific sections of the network during peak hours of a day. At the same time, since the restriction may significantly increase the truck company's cost, the MPO could compensate the truck company in order to smooth the implementation of the restriction.

Some interesting extensions of this study include but not limited to the following: (1) comprehensively investigate the dynamic Braess-like paradox based on the proposed modeling framework and discuss the necessary and sufficient conditions for the paradox to occur; (2) robustify the problem by incorporating uncertain travel demand; (3) model the upper level problem as service network design problem which is more realistic yet more challenging; (4) assume that there exist multiple truck companies competing with each other and study the resulting equilibrium problem with equilibrium constraints (EPEC).

**ACKNOWLEDGMENT**

This research was supported in part by the grant award CMMI-0900040 from the National Science Foundation (NSF) and the grant awards from the Mid-Atlantic Universities Transportation Center (MAUTC). The work described in this paper has been presented in Institute for Operations Research and the Management Sciences (INFORMS) Transportation Science and Logistics (TSL) Society Workshop 2011 and the 4th International Symposium on Dynamic Traffic Assignment (DTA 2012). The authors greatly acknowledge the helpful and constructive comments from the audience at these two conferences.